\documentclass[leqno, nowrap]{amsart}

\usepackage{amsmath}
\numberwithin{equation}{section}

\usepackage{amsfonts}
\usepackage{amssymb}
\usepackage[top=1.5in, bottom=1.5in, left=1.5in, right=1.5in]{geometry}

\usepackage{amsthm}
\newtheorem{Theorem}{Theorem}[section]
\newtheorem{Lemma}[Theorem]{Lemma}
\newtheorem{Corollary}[Theorem]{Corollary}
\newtheorem{Proposition}[Theorem]{Proposition}
\newtheorem*{Remark}{Remark}

\usepackage{indentfirst}
\usepackage{verbatim}
\usepackage{endnotes}

\newcommand{\cal}{\mathcal}
\newcommand{\R}{\mathbb{R}}

\renewcommand{\H}{\cal{H}}

\newcommand{\del}{\partial}

\newcommand{\grad}{\nabla}
\renewcommand{\div}{\mathrm{div}}
\newcommand{\beq}{\begin{equation}}
\newcommand{\eeq}{\end{equation}}
\newcommand{\bTh}{\begin{Theorem}}
\newcommand{\eTh}{\end{Theorem}}
\newcommand{\bCo}{\begin{Corollary}}
\newcommand{\eCo}{\end{Corollary}}
\newcommand{\bLem}{\begin{Lemma}}
\newcommand{\eLem}{\end{Lemma}}
\newcommand{\bProp}{\begin{Proposition}}
\newcommand{\eProp}{\end{Proposition}}
\newcommand{\bProof}{\begin{proof}}
\newcommand{\eProof}{\end{proof}}
\newcommand{\bRe}{\begin{Remark}}
\newcommand{\eRe}{\end{Remark}}
\newcommand{\Hs}[2]{{H^{#1}(#2)}}

\newcommand{\Lp}[2]{{L^{#1}(#2)}}

\newcommand{\leC}{\lesssim}

\newcommand{\eps}{\varepsilon}
\renewcommand{\phi}{\varphi}
\newcommand{\lam}{\lambda}

\renewcommand{\and}{\text{and}}

\newcommand{\sgn}{\mathrm{sgn}}

\begin{document}

\hyphenpenalty=1000000

\title[Nodal sets of Dirichlet and Neumann eigenfunctions]{Lower bounds for nodal sets of Dirichlet and Neumann eigenfunctions}
\author{Sinan Ariturk}
\begin{abstract}
Let $\phi$ be a Dirichlet or Neumann eigenfunction of the Laplace-Beltrami operator on a compact Riemannian manifold with boundary.
We prove lower bounds for the size of the nodal set $\{ \phi = 0 \}$.
\end{abstract}

\maketitle

\section{Introduction}

Let $(M,g)$ be a compact smooth Riemannian manifold with boundary.
Let $\Delta$ be the Laplace-Beltrami operator.
Let $\lam \ge 1$ and let $\phi$ be an eigenfunction of $-\Delta$, i.e. a smooth real-valued function on $M$ with
\[
	-\Delta \phi = \lam \phi
\]
over the interior of $M$.
We will assume that  $\phi$ is a Dirichlet eigenfunction, meaning
\[
	\phi \Big|_{\del M} = 0
\]
or a Neumann eigenfunction, meaning
\[
	\del_\nu \phi \Big|_{\del M} = 0
\]
where $\nu$ is the outward unit normal vector on $\del M$ and $\del_\nu$ is the corresponding directional derivative.
Define the nodal set
\[
	Z = \bigg\{  x \in M : \phi(x)=0, x \notin \del M \bigg\}
\]
Let $n$ be the dimension of $M$ and let $\H$ be the $(n-1)$-dimensional Hausdorff measure on $M$.
We will prove lower bounds for $\H(Z)$.

We use the notation $A \leC B$ to mean there is a positive constant $C$, independent of $\lam$ and $\phi$, such that $A \le CB$.

\bTh
\label{th}
If $\phi$ is a Neumann eigenfunction, then
\[
	\lam^\frac{5-2n}{6} \leC \H(Z)
\]
If $\phi$ is a Dirichlet eigenfunction and $n \le 3$, then
\[
	\lam^\frac{5-2n}{6} \leC \H(Z)
\]
If the boundary is strictly geodesically concave and $\phi$ is a Dirichlet eigenfunction, then for $n \le 4$,
\[
	\lam^\frac{3-n}{4} \leC \H(Z)
\]
\eTh

If $(M,g)$ is a compact real analytic Riemannian manifold with boundary, then Donnelly and Fefferman \cite{DF} proved that
\[
	\lam^{1/2} \leC \H(Z) \leC \lam^{1/2}
\]
If $(M,g)$ is a compact smooth Riemannian manifold without boundary, then Colding and Minicozzi \cite{CM} proved that
\beq
\label{CMin}
	\lam^\frac{3-n}{4} \leC \H(Z)
\eeq
This same result was later obtained by Hezari and Sogge \cite{HS}.
Their argument was based on the identity
\beq
\label{SZid}
	\lam \int_M | \phi | \,dV = 2 \int_Z | \grad \phi | \,dS
\eeq
where $dV$ is the Riemannian volume measure and $dS$ is the Riemannian surface measure on $Z$.
This identity had been proven by Sogge and Zelditch \cite{SZ}, who also showed that
\beq
\label{SZin}
	\lam^{-\frac{n-1}{8}} \leC \int_M | \phi | \,dV
\eeq
Hezari and Sogge \cite{HS} proved that
\beq
\label{HSin}
	\int_Z | \grad \phi |^2 \,dS \leC \lam^{3/2}
\eeq
and then used \eqref{SZid}, \eqref{SZin}, and \eqref{HSin} to obtain the bound \eqref{CMin}.

We will prove analogues of \eqref{SZid}, \eqref{SZin}, and \eqref{HSin} for a compact smooth Riemannian manifold with boundary.
This will enable us to establish Theorem \ref{th}.
In particular, we will prove the following.

\bTh
\label{SZ}
If $\phi$ is a Dirichlet or Neumann eigenfunction, then
\[
	\lam \int_M | \phi | \,dV = \int_{\del M} | \del_\nu \phi | \,dS + 2 \int_Z | \grad \phi | \,dS
\]
More generally, for any function $f$ in $C^2(M)$,
\[
	\int_M \Big( (\Delta + \lam)f \Big) | \phi | \,dV = \int_{\del M} f | \del_\nu \phi | \,dS + \int_{\del M} |\phi| \del_\nu f \,dS + 2 \int_Z f | \grad \phi | \,dS
\]
\eTh

For a Neumann eigenfunction, the first term on the right side is zero, and this identity is the same as \eqref{SZid}.
For a Dirichlet eigenfunction, the integral over $\del M$ is an additional obstacle and causes the argument to break down in higher dimensions.

\subsection*{Acknowledgements}
	I would like to thank Christopher Sogge for suggesting this problem and for his invaluable guidance.

\section{Proofs}

Define
\[
	P = \bigg\{ x \in M : \phi(x) > 0, x \notin \del M \bigg\}
\]
and
\[
	N = \bigg\{ x \in M : \phi(x) < 0, x \notin \del M \bigg\}
\]
We can write $M$ as a disjoint union
\[
	M = P \cup N \cup \del M \cup Z
\]
Define
\[
	\Omega = \bigg\{ x \in M: \phi(x)=0 \bigg\}
\]
and
\[
	\Sigma = \bigg\{ x \in \Omega : \grad \phi(x)=0 \bigg\}
\]

\bLem
If $\phi$ is a Dirichlet or Neumann eigenfunction, then $\H(\Omega) < \infty$, and the Hausdorff dimension of $\Sigma$ is at most $n-2$.
If $\phi$ is a Neumann eigenfunction, then the Hausdorff dimension of $\Omega \cap \del M$ is at most $n-2$.
\eLem

\bProof
Fix a point $p$ in $M$.
To prove the first statement, it suffices to find a neighborhood $U$ of $p$ in $M$ such that $\H(\Omega \cap U) < \infty$ and $\Sigma \cap U$ has Hausdorff dimension at most $n-2$.
If $\phi(p) \neq 0$, then finding such a neighborhood $U$ is trivial.
So we assume $\phi(p)=0$.
By Donnelly and Fefferman \cite{DF}, the eigenfunction $\phi$ only vanishes to finite order at $p$.
If $p$ is in the interior of $M$, we use geodesic normal coordinates about $p$.
Then by Hardt and Simon~\cite{HaSi}, we can obtain $U$.

If $p$ is on the boundary $\del M$, then we use boundary normal coordinates $(x_1,\ldots, x_n)$ about $p$.
These are defined by first letting $(x_1, \ldots x_{n-1})$ be geodesic normal coordinates on $\del M$ about $p$, with respect to the metric on $\del M$ induced by $g$.
Then for fixed $x_1, \ldots, x_{n-1}$, the curves $x_n \to (x_1,\ldots, x_n)$, for $x_n \ge 0$, are geodesics in $M$ which intersect $\del M$ normally.
These coordinates are well-defined near $p$ and allow us to identify some neighborhood of $p$ with
\[
	B_+ = \bigg\{ x \in \R^n : |x| < \eps, x_n \ge 0 \bigg\}
\]
for some small $\eps>0$.
Here the point $p$ is being identified with the origin in $\R^n$. 
Let $g_{ij}$ be the Riemannian metric on $B_+$.
Let
\[
	B = \bigg\{ x \in \R^n : |x| < \eps \bigg\}
\]
We extend the metric $g_{ij}$ to $B$ so that it is even in the $x_n$-variable.
Let $g^{ij}$ be the cometric, defined so that the matrix $[g^{ij}]$ is the inverse matrix of $[g_{ij}]$.
Define
\[
	J = \Big( \det [g_{ij}] \Big)^{1/2}
\]
The functions $g_{ij}$, $g^{ij}$, and $J$ are Lipschitz continuous and bounded on $B$.
If $\phi$ is a Dirichlet eigenfunction, extend $\phi$ to $B$ so that it is odd in the $x_n$-variable.
If $\phi$ is a Neumann eigenfunction, extend $\phi$ to $B$ so that it is even in the $x_n$-variable.
Then the extended function $\phi$ is in $C^1(B) \cap H^2(B)$.
Let $\psi$ be a smooth function on $\R^2$ with compact support contained strictly inside $B$.
By Green's identity,
\[
	\sum_{i,j=1}^n \int_B (D_j \phi) (D_i \psi) J g^{ij} \,dx = \int_B \lam \phi \psi J \,dx
\]
That is,
\[
	\bigg( \sum_{i,j=1}^n D_i J g^{ij} D_j \phi \bigg) + \lam J \phi = 0
\]
We can write this equation as
\[
	\bigg( \sum_{i,j=1}^n J g^{ij} D_i D_j \phi + (D_i J g^{ij}) D_j \phi \bigg) + \lam J \phi = 0
\]
Now by Hardt and Simon \cite{HaSi}, we can obtain $U$.

It remains to prove the second statement.
Fix a point $p$ in $(\Omega \setminus \Sigma) \cap \del M$.
It suffices to show that there is a neighborhood $V$ of $p$ in $\del M$ such that the Hausdorff dimension of $(\Omega \setminus \Sigma) \cap V$ is at most $n-2$.
The set $\Omega \setminus \Sigma$ is a hypersurface with normal vector $\grad \phi(p)$ at $p$.
Since $\phi$ is a Neumann eigenfunction and $\grad \phi(p) \neq 0$, the sets $\Omega \setminus \Sigma$ and $\del M$ intersect transversally, which yields $V$.
\eProof

In particular, it follows that $\del P$ is smooth almost everywhere, with respect to $\H$, so the divergence theorem and Green's identities hold on $P$.
See, e.g., Evans and Gariepy \cite{EG}.
Let $\eta$ be the outward unit normal on $\del P$, defined at these smooth points, and let $\del_\eta$ be the corresponding directional derivative.
On $Z \setminus \Sigma$, we have
\[
	\eta = -\frac{\grad \phi}{| \grad \phi|}
\]
At any point on $\del M \cap \del P$ where $\eta$ is defined, we have
\[
	\eta = \nu
\]

\bProof[Proof of Theorem \ref{SZ}]
By Green's identity,
\[
\begin{split}
	\int_P \Big( (\Delta + \lam)f \Big) | \phi | \,dV &= \int_P \Big( (\Delta + \lam)f \Big) \phi \,dV
	\\
	&= \int_P f (\Delta + \lam) \phi  \,dV - \int_{\del P} f \del_\eta \phi \,dS + \int_{\del P} \phi \del_\eta f \,dS
	\\
	&= - \int_{\del P \cap \del M} f \del_\eta \phi \,dS - \int_Z f \del_\eta \phi \,dS + \int_{\del P \cap \del M} \phi \del_\eta f \,dS
	\\
	&= \int_{\del P \cap \del M} f |\del_\nu \phi| \,dS + \int_Z f | \grad \phi | \,dS + \int_{\del P \cap \del M} |\phi| \del_\nu f \,dS
\end{split}
\]
The last equality holds because $-\del_\eta \phi = | \del_\nu \phi |$ over $\del P \cap \del M$ and $-\del_\eta \phi = | \grad \phi |$ over $\del P \cap Z$.
We can similarly obtain
\[
	\int_N \Big( (\Delta + \lam)f \Big) | \phi | \,dV = \int_{\del N \cap \del M} f |\del_\nu \phi| \,dS + \int_Z f | \grad \phi | \,dS + \int_{\del N \cap \del M} |\phi| \del_\nu f \,dS
\]
Now
\[
\begin{split}
	\int_M \Big( (\Delta + \lam)f \Big) | \phi | \,dV
		&= \int_P \Big( (\Delta + \lam)f \Big) | \phi | \,dV + \int_N \Big( (\Delta + \lam)f \Big) | \phi | \,dV
	\\
	&= \int_{\del M} f | \del_\nu \phi | \,dS + \int_{\del M} |\phi| \del_\nu f \,dS + 2 \int_Z f | \grad \phi | \,dS
\end{split}
\]
\eProof

The following lemma is an analogue of \eqref{SZin}.
\bLem
\label{L1}
If $\phi$ is a Dirichlet or a Neumann eigenfunction, then
\[
	\lam^\frac{1-n}{6} \leC \| \phi \|_\Lp{1}{M}
\]
If the boundary is strictly geodesically concave and $\phi$ is a Dirichlet eigenfunction, then
\[
	\lam^\frac{1-n}{8} \leC \| \phi \|_\Lp{1}{M}
\]
\eLem

\bProof
Fix $p$ satisfying $2<p<\frac{2(n+1)}{n-1}$.
Then, by Smith \cite{S},
\beq
\label{Lp}
	\| \phi \|_\Lp{p}{M} \leC \lam^\frac{(n-1)(p-2)}{6p}
\eeq
If the boundary is strictly geodesically concave and $\phi$ is a Dirichlet eigenfunction, then by Grieser \cite{G} and Smith-Sogge \cite{SS},
\[
	\| \phi \|_\Lp{p}{M} \leC \lam^\frac{(n-1)(p-2)}{8p}
\]
Let $\theta = \frac{p-2}{2(p-1)}$.
By H\"older's inequality,
\[
	1 = \| \phi \|_\Lp{2}{M} \le \| \phi \|_\Lp{1}{M}^\theta \| \phi \|_\Lp{p}{M}^{1-\theta}
\]
The estimates now follow.
\eProof

\bRe
On the flat unit disc $\{ |x| \le 1 \}$ in $\R^2$, there are whispering gallery modes, which are concentrated near the boundary.
It follows from Grieser \cite{G} that Lemma \ref{L1} is sharp for these eigenfunctions.
However, for $n \ge 3$, Smith and Sogge \cite{SSbd} conjectured that \eqref{Lp} can be strengthened to
\beq
\label{Lpconj}
	\| \phi \|_\Lp{p}{M} \leC \lam^\frac{(3n-2)(p-2)}{24p}
\eeq
Applying H\"older's inequality as above would then yield
\[
	\lam^\frac{2-3n}{24} \leC \| \phi \|_\Lp{1}{M}
\]
\eRe

The following lemma is an analogue of \eqref{HSin}.

\bLem
\label{HSinbd}
If $\phi$ is a Dirichlet or Neumann eigenfunction, then
\[
	\int_Z | \grad \phi |^2 \,dS \leC \lam^{3/2}
\]
\eLem

\bProof
This will follow from the identity
\[
	-\int_M \sgn(\phi) \, \div \Big(| \grad \phi | \grad \phi \Big) \,dV = \int_{\del M} | \del_\nu \phi |^2 \,dS + 2 \int_Z | \grad \phi |^2 \,dS
\]
We first prove this identity.
Note that $-\del_\eta \phi = | \grad \phi |$ over $Z \setminus \Sigma$.
If $\phi$ is a Dirichlet eigenfunction, then we also have $| \grad \phi | = -\del_\eta \phi = | \del_\nu \phi |$ at any point on $\del P \cap \del M$ where $\eta$ is defined.
By the divergence theorem,
\[
\begin{split}
	- \int_P \div \Big(| \grad \phi | \grad \phi \Big) \,dV 
	&= - \int_{\del P} | \grad \phi | \del_\eta \phi \,dS
	\\
	&=\int_{\del P \cap \del M} | \del_\nu \phi |^2 \,dS + \int_{Z} | \grad \phi |^2 \,dS
\end{split}
\]
Similarly,
\[
	\int_N \div \Big(| \grad \phi | \grad \phi \Big) \,dV = \int_{\del N \cap \del M} | \del_\nu \phi |^2 \,dS + \int_{Z} | \grad \phi |^2 \,dS
\]
Adding these equations establishes the identity.
Now we have
\[
\begin{split}
	\int_Z | \grad \phi |^2 \,dS &\le \int_M \Big| \div \Big(| \grad \phi | \grad \phi \Big) \Big| \,dV
		\\&\leC \| \phi \|_\Hs{2}{M} \| \phi \|_\Hs{1}{M}
		\\&\leC \lam^{3/2}
\end{split}
\]
\eProof

For a Dirichlet eigenfunction, we also need the following lemma.
\bLem
\label{deltrace}
If $\phi$ is a Dirichlet eigenfunction, then
\[
	\bigg( \int_{\del M} | \del_\nu \phi |^2 \,dS \bigg)^{1/2} \leC \lam^{1/2}
\]
\eLem

This lemma follows from a much more general result obtained by Tataru \cite{DT}.
There is also the following short proof.

\bProof
Let $X$ be a smooth first-order differential operator on $M$ with $X=\del_\nu$ over $\del M$.
Then, by Green's identity,
\[
\begin{split}
	\int_M u [X, \Delta] u \,dV &= - \lam \int_M u Xu \,dV - \int_M u \Delta X u \,dV
		\\&= \int_M (\Delta u) (Xu) \,dV - \int_M u \Delta X u \,dV
		\\&= \int_{\del M} (\del_\nu u) (Xu) \,dS
		\\&= \int_{\del M} | \del_\nu u |^2 \,dS
\end{split}
\]
Since $[X,\Delta]$ is a second-order differential operator, this yields
\[
\begin{split}
	\int_{\del M} | \del_\nu u |^2 \,dS &= \int_M u [X, \Delta] u \,dV
		\\&\leC \| u \|_\Lp{2}{M} \| u \|_\Hs{2}{M}
		\\&\leC \lam
\end{split}
\]
\eProof

We can now prove Theorem \ref{th}.

\bProof[Proof of Theorem \ref{th}]
First assume $\phi$ is a Neumann eigenfunction.
By Theorem \ref{SZ} and Lemma \ref{HSinbd},
\[
	\lam \int_M | \phi | \,dV = 2 \int_Z | \grad \phi | \,dS \leC \H(Z)^{1/2} \lam^{3/4}
\]
We can rewrite this as
\[
	\lam^{1/2} \bigg( \int_M | \phi | \,dV \bigg)^2 \leC \H(Z)
\]
So by Lemma \ref{L1},
\[
	\lam^\frac{5-2n}{6} \leC \H(Z)
\]

Now assume $\phi$ is a Dirichlet eigenfunction.
By Theorem \ref{SZ}, Lemma \ref{HSinbd}, and Lemma~\ref{deltrace},
\[
	\lam \int_M | \phi | \,dV = \int_{\del M} | \del_\nu \phi | \,dS + 2 \int_Z | \grad \phi | \,dS \leC \lam^{1/2} + \H(Z)^{1/2} \lam^{3/4}
\]
We can rewrite this as
\[
	\lam^{1/2} \bigg( \int_M | \phi | \,dV \bigg)^2 \leC \H(Z) + \lam^{-1/2}
\]
Now applying Lemma \ref{L1} yields the desired estimates.
\eProof

\bRe
If \eqref{Lpconj} is true, then we would have a better lower bound for the $L^1$ norm of $\phi$.
If $\phi$ is a Neumann eigenfunction, this would yield
\[
	\lam^\frac{8-3n}{12} \leC \H(Z)
\]
The same bound would hold if $\phi$ is a Dirichlet eigenfunction and $n \le 4$.
\eRe

\linespread{1.3}


\begin{thebibliography}{9}

\bibitem{CM}T. H. Colding and W. P. Minicozzi, II. \emph{Lower bounds for nodal sets of eigenfunctions}, arXiv:1009.4156, to appear in Comm. Math. Phys.

\bibitem{DF}H. Donnelly and C. Fefferman. \emph{Nodal sets of eigenfunctions: Riemannian manifolds with boundary}. Analysis, et cetera, 251-262, Academic Press, Boston, MA, 1990.

\bibitem{EG}L. C. Evans and R. F. Gariepy. \emph{Measure theory and fine properties of functions}. Studies in Advanced Mathematics. CRC Press, Boca Raton, FL, 1992.

\bibitem{G}D. Grieser. \emph{$L^p$ bounds for eigenfunctions and spectral projections of the Laplacian near concave boundaries}. Ph. D. Thesis. University of California, Los Angeles: USA, 1992.

\bibitem{HaSi}R. Hardt and L. Simon. \emph{Nodal sets for solutions of elliptic equations}. J. Differential Geom. 30 (1989), no. 2, 505-522.

\bibitem{HS}H. Hezari and C. D. Sogge. \emph{A natural lower bound for the size of nodal sets}, arXiv:1107.3440, preprint.

\bibitem{S}H. F. Smith. \emph{Sharp $L^2 \to L^q$ bounds on spectral projectors for low regularity metrics}. Math. Res. Lett. 13 (2006), no 5-6, 967-974.

\bibitem{SS}H. F. Smith and C. D. Sogge. \emph{On the critical semilinear wave equation outside convex obstacles}. J. Amer. Math. Soc. 8 (1995), no. 4, 879-916.

\bibitem{SSbd}H. F. Smith and C. D. Sogge. \emph{On the $L^p$ norm of spectral clusters for compact manifolds with boundary}. Acta Math. 198 (2007), 107-153.

\bibitem{SZ}C. D. Sogge and S. Zelditch. \emph{Lower bounds on the Hausdorff measure of nodal sets}, Math. Res. Lett. 18 (2011), no. 1, 25-37.

\bibitem{DT}D. Tataru. \emph{On the regularity of boundary traces for the wave equation}. Ann. Scuola Norm. Sup. Pisa Cl. Sci. (4) 26 (1998), no. 1, 185-206.

\end{thebibliography}
\end{document}